\documentclass[11pt]{amsart}
\usepackage{}
\usepackage{amssymb}
\theoremstyle{plain}
\newtheorem{Thm}{Theorem}

\newtheorem{Pro}[Thm]{Proposition}

\begin{document}

\title[gap theorems]
{gap theorems for locally conformally flat manifolds}

\author{Li Ma}

\address{ Dr.Li Ma, Distinguished professor\\
Department of mathematics \\
Henan Normal university \\
Xinxiang, 453007 \\
China}

\email{nuslma@gmail.com}

\thanks{ The research is partially supported by the National Natural Science
Foundation of China 11271111 and SRFDP 20090002110019}

\begin{abstract}
In this paper, we prove a gap result for a locally conformally flat
complete non-compact Riemannian manifold with bounded non-negative
Ricci curvature and a scalar curvature average condition. We show
that if it has positive Green function, then it is flat. This result
is proved by setting up new global Yamabe flow. Other extensions
related to bounded positive solutions to a schrodinger equation are
also discussed. A global existence of Yamabe flow on a locally
conformally flat complete non-compact Riemannian manifold with
bounded non-negative sectional curvature is proved.

{ \textbf{Keywords}: global Yamabe flow, Poisson equation, Harnack
inequality, gap theorem}

 { \textbf{Mathematics Subject
Classification 2000}: 35J60, 53J60,58J05}
\end{abstract}

 \maketitle

\section{Introduction} In the recent work of L.Ni
\cite{ni}, the author proves a very interesting gap theorem on a
complete noncompact Kaehler manifold with nonnegative holomorphic
bisectional curvature via the use heat equation technique. The
remarkable feature in Ni's result is the fast quadratic decay
condition about the average of the scalar curvature in large balls.
His result can be considered as some kind of  positive mass theorem.
In this paper, we consider other scalar curvature conditions and we
prove the following new gap theorem on a locally conformally flat
(LCF in short) complete non-compact Riemannian manifold with bounded
non-negative Ricci curvature.

\begin{Thm}\label{main} Assume that $(M^n,g_0)$ ($n\geq 3$) is a
locally conformally flat complete non-compact Riemannian manifold
with bounded non-negative Ricci curvature and with the scalar
curvature average condition for some $x_0\in M^n$,
\begin{equation}\label{average}
\int_0^\infty\frac{r}{Vol_{x_0}(r)}\int_{B_{x_0}(r)}R_0(y)dv_{g_0}dr<\infty,
\end{equation}
where $R_0(y)$ is the scalar curvature function of $(M^n,g_0)$ and
$Vol_{x_0}(r)$ is the volume of the ball $B(x_0,r)$ in $(M^n,g_0)$.
Assume further that it is non-parabolic, i.e., it has a positive
Green function. Then it is flat.
\end{Thm}

We now give two remarks about the proof and the conditions of
Theorem \ref{main}.

\textbf{Remark 1}. According to heat kernel estimate of Li-Yau
\cite{ly}(1986) and the work of N. Varopoulos \cite{V}, we know that
for the complete non-compact Riemannian manifold with non-negative
Ricci curvature, the non-parabolic condition is equivalent to the
ball volume assumption that
\begin{equation}\label{volume-control}
\int_1^\infty\frac{r}{Vol_{x_0}(r)}dr<\infty,
\end{equation}
where $x_0\in M:=M^n$ is some point. We remark that on the
non-parabolic manifold $(M,g_0)$, if the scalar curvature average
condition at some point $x_0$, then it is also true at any point $x$
of $M$. Further extensions of our result above are stated in the
last section.

\textbf{Remark 2}. The key part of our proof of Theorem \ref{main}
is our discovery that there is a remarkable relation between the
Poisson equation (or Schrodinger equation) and the existence of the
global Yamabe flow on complete non-compact Riemannian manifold with
non-negative scalar curvature.

The idea of the proof of Theorem \ref{main} above is below. Assume
that $(M,g_0)$ is non-flat. Then the scalar curvature $R_0(x)$ of it
is non-trivial and non-negative. We use the non-parabolic condition
and the average condition \ref{average} to solve the Poisson
equation
$$
-\Delta w(x)=R_0(x), \ \ in \ \ M
$$
and get a non-negative solution $w$. Then we use the function $w$ to
get a lower barrier of the Yamabe flow to prevent it from locally
collapsing and show that we have a global Yamabe flow. We then use
the fundamental Harnack inequality of B.Chow \cite{chow} for locally
conformally flat manifolds to show that the scalar curvature must be
trivial. This then concludes the result of Theorem \ref{main}.

From our argument of Theorem \ref{main}, we can actually prove
stronger results. See Theorems \ref{main+1} and \ref{main+2} in the
last section.

The result above can be considered as a generalized positive mass
theorem \cite{SY} in the sense that if $(M,g_0)$ is a locally
conformally flat complete non-compact Riemannian manifold with
bounded non-trivial non-negative Ricci curvature and the condition
(\ref{average}), then it is parabolic, i.e., the Green function must
be negative somewhere. It is reasonable to believe that our result
above should have more applications to locally conformally flat
geometric structures. In our previous work \cite{ma}, we can use the
positive mass theorem to solve an open question posed in
\cite{chow06}.

There are many gap theorems for complete non-compact Riemannian
manifold with suitable assumption about curvatures stronger than
non-negative Ricci curvature. For more references, one may see
L.Ni's paper \cite{ni}. See also \cite{chow06} and \cite{shi}. Here
we mention the following result of Ma-Cheng \cite{macheng} as a
comparison.

\begin{Thm}\label{ma} Assume that $(M,g_0)$ is a
locally conformally flat complete non-compact Riemannian manifold
with bounded non-negative Ricci curvature. Assume further that it
has a Ricci pinching condition in the sense that there is a positive
constant $\epsilon$ such that $$ Rc(g_0)(x)\geq \epsilon
R_0(x)g_0(x).
$$ Then it is flat.
\end{Thm}

Recall that the asymptotic volume ratio (AVR) of a complete
noncompact Riemannian manifold $(M^n,g)$ with non-negative Ricci
curvature is defined by
$$
AVR(g)=\lim_{r\to\infty}\frac{Vol(B_r(0))}{r^n}.
$$
Note that the definition does not depends on the choice of the base
point $0\in M$.

Then we can show the following global Yamabe flow result.

\begin{Thm}\label{mai}
Assume that $(M,g_0)$ is a locally conformally flat complete
non-compact Riemannian manifold with bounded non-negative sectional
curvature and its curvature tends to zero at infinity. Assume
further that $AVR(g_0)>0$. Then there is a global Yamabe flow $g(t)$
on $(M,g_0)$ with $g_0$ as the initial metric.
\end{Thm}
We remark that it is generally conjectured that there is a global
Yamabe flow $g(t)$ on $(M,g_0)$ locally conformally flat complete
non-compact Riemannian manifold with bounded non-negative Ricci
curvature and its curvature tends to zero at infinity. This is still
an open problem, at least for the locally conformally flat
$M=S^{n-1}\times R$. For further interesting classification result,
one may see \cite{DS}.

The plan of this paper is below. In section \ref{sect2}, we discuss
the existence of a non-negative weak solution to the Poisson
equation in $(M,g_0)$ under the condition \ref{volume-control}. In
section \ref{sect3}, we introduce the global Yamabe flow with the
initial data $g_0$. In section \ref{sect4}, we prove the main result
and discuss other extensions. In the last section we prove Theorem
\ref{mai}. We may use $C$ to denote various uniform constants.

\section{Poisson equation}\label{sect2}

Let $(M,g_0)$ be a complete non-compact Riemannian manifold of
dimension $n\geq 3$ and with non-negative Ricci curvature and the
scalar curvature average condition (\ref{average}).

We assume (\ref{volume-control}). Then according to Li-Yau \cite{ly}
(1986) we know that $(M,g_0)$ is non-parabolic, i.e., there is a
positive Green function $G(x,y)$ such that
$$
C^{-1}\int_{d(x,y)}^\infty\frac{r}{Vol_{x_0}(r)}dr\leq G(x,y)\leq
C\int_{d(x,y)}^\infty\frac{r}{Vol_{x_0}(r)}dr
$$
for some uniform dimension constant $C>0$, where $d(x,y)$ is the
distance function between the points $x$ and $y$ in $(M,g_0)$.

Given any non-negative continuous function $f:M\to R$ such that for
any $x\in M$,
$$
\int_0^\infty\frac{r}{Vol_{x}(r)}\int_{B_x(r)}f(y)dv_{g_0}dr<\infty.
$$
Then the non-negative function
$$
w(x)=\int_MG(x,y)f(y)dv_{g_0}
$$
is well-defined and is the non-negative solution in the class
$C^{1,\alpha}_{loc}(M)\cap W^{2,p}_{loc}(M)$ ($\alpha\in (0,1)$ and
$1<p<\infty$)to the Poisson equation
$$
-\Delta_{g_0}w=f, \ \ in \ \ M.
$$

In fact, using the Fubini theorem, we know that
$$
\int_MG(x,y)f(y)dv_{g_0}=\int_0^\infty dr \int_{\partial
B(x,r)}G(x,y)f(y)d\sigma_y,
$$
which is no bigger than
$$
C\int_0^\infty dr\int_{r}^\infty\frac{s}{Vol_{x}(s)}ds\int_{\partial
B(x,r)}f(y)d\sigma_y,
$$
where $C$ is some uniform constant. Exchanging the order of $s$ and
$r$, the latter term can be written as
$$
C\int_0^\infty ds\int_{0}^s\frac{s}{Vol_{x}(s)}\int_{\partial
B(x,r)}f(y)d\sigma_ydr,
$$
which is
$$
C\int_0^\infty\frac{s}{Vol_{x}(s)}\int_{B_x(s)}f(y)dv_{g_0}ds<\infty.
$$

Since for any $\phi\in C_0^2(M)$, we have
$$
\int_M\Delta \phi(x)\int_MG(x,y)f(y)dv_{g_0}=\int_Mf(y)
\int_MG(x,y)\Delta\phi(x)dv_{g_0},
$$
and the latter is
$$
\int_Mf(y)\int_M\Delta G(x,y)\phi(x)dv_{g_0}=-\int_Mf(y)\phi(y)
$$ after a use of integration by part. Hence we have
$$
-\Delta w=f
$$
holds weakly. Using the Schauder theory \cite{chow06} we know that
$w\in C^{1,\alpha}_{loc}(M)$ for any $\alpha\in (0,1)$. Using the
Calderon-Zugmund $L^p$ theory \cite{m} we know that $w\in
W^{2,p}_{loc}(M)$. Higher regularity of $w$ can be obtained by
assumption of higher regularity of the function $f$.

In summary we have proved the following result, which is well known
to experts in geometric analysis.
\begin{Pro}\label{ni} Let $(M,g_0)$ be a complete non-compact Riemannian manifold of
dimension $n\geq 3$ and with non-negative Ricci curvature and the
assumption (\ref{volume-control}). Then for any non-negative
continuous function $f:M\to R$ such that for any $x\in M$,
$$
\int_0^\infty\frac{r}{Vol_{x}(r)}\int_{B_x(r)}f(y)dv_{g_0}dr<\infty,
$$
there is the non-negative solution $w$ in the class
$C^{1,\alpha}_{loc}(M)\cap W^{2,p}_{loc}(M)$ ($\alpha\in (0,1)$ and
$1<p<\infty$)to the Poisson equation $ -\Delta_{g_0}w=f $ in $M$.
Furthermore, if $f\in C^\infty$, then $w\in C^\infty$.
\end{Pro}

\section{global Yamabe flow with positive scalar curvature}\label{sect3}
In this section, we assume that $(M,g_0)$ is a complete non-compact
Riemannian manifold with non-negative scalar curvature $R_0$. Assume
that the Poisson equation $-\Delta w=\frac{n-2}{4(n-1)}R_0$ on $M$
has a non-negative solution $w$. We emphasis that in this section,
we don't need the conditions that the scalar curvature is bounded
and the manifold is locally conformally flat. See also \cite{AM}.

We may assume that $(M,g_0)$ is given with $R(g_0)\geq 0$ being
non-trivial.

 Recall that the
Yamabe flow is a family of pointwise conformally equivalent metrics
$g(t)$ ($t\in [0,T)$, $T>0$) evolved by the evolution equation
\begin{equation}\label{yamabe0}
\partial_t g(t)=-R(t)g(t), \ \ on \ \ M\times (0,T),
\end{equation}
with the initial metric $g(0)=g_0$. Here $R(t)=R(g(t))$ is the
scalar curvature of the metric $g(t)$. We shall write $R_0(x)$ by
the scalar curvature of $g_0$. Under the Yamabe flow, we know that
\begin{equation}\label{scalar}
\partial_tR(t)=(n-1)\Delta_{g(t)}R+R^2, \ in \ M\times (0,T).
\end{equation}

Let $g(t)=u(x,t)^{4/(n-2)}g_0$ for positive smooth functions
$u(x,t)$ and let $p=\frac{n+2}{n-2}$. Then we have that
\begin{equation}\label{rule}
R(t)=u^{-p}[-\frac{4(n-1)}{n-2}\Delta u+R_0u].
\end{equation}

Then the evolution equation (\ref{yamabe0}) becomes the parabolic
equation for the positive functions $u(x,t)$
\begin{equation}\label{yamabe}
\partial_tu^p=(n-1)p[\Delta u-\frac{n-2}{4(n-1)}R_0(x)u], \ on \ M\times (0,T),
\end{equation}
with the initial value $u(x,0)=1$. Here $\Delta$ is the Laplace
operator of the metric $g_0$. We shall write by $\nabla u$ the
gradient of the function $u$ in the metric $g_0$.

We know from the standard parabolic theory that for any bounded
smooth domain $\Omega\subset M$, there exists some small $T>0$ such
that the equation (\ref{yamabe}) on $\Omega\times (0,T)$ always has
a positive solution $\Omega\times (0,T)$ with the initial and
boundary conditions $u(x,t)=1$ either for $t>0, x\in \partial\Omega$
or $t=0, x\in \Omega$. Since $R_0\geq 0$, we have
$$
\partial_tu^p\leq (n-1)p\Delta u, \ in \ \Omega\times (0,T)
$$
which implies that $u(x,t)\leq 1$ in $\Omega\times [0,T)$.

With the abuse notation, we still use $g(t)=u(x,t)^{4/(n-2)}g_0$ on
$\Omega\times (0,T)$ for the locally defined solution $u$.

Since $R_0\geq 0$ on $M$, by the equation (\ref{scalar}) on
$\Omega\times (0,T)$ and the maximum principle \cite{m} we know that
$R(t)>0$ on $\Omega\times (0,T)$. Fix any $t\in (0,T)$. Using the
formula (\ref{rule}) we have
\begin{equation}\label{schrodinger}
-\frac{4(n-1)}{n-2}\Delta u+R_0u>0, \ in \Omega.
\end{equation}
We may rewrite this as
$$
\frac{n-2}{4(n-1)}R_0>\frac{\Delta u}{u}.
$$

 Set $v=\log u$. Then we have
$$
\Delta v=\frac{\Delta u}{u}-\frac{|\nabla u|^2}{u^2}\leq
\frac{\Delta u}{u}.
$$
Hence, we have
$$
\Delta v<\frac{n-2}{4(n-1)}R_0, \ in \ \Omega.
$$

We now let $f(x)=\frac{n-2}{4(n-1)}R_0$ solve the Poisson equation
$\Delta (-w)=f$ for $w$ in the whole manifold $M$. Then we have
$$
\Delta (v+w)<0, \ in \ \Omega
$$
with the boundary condition $v+w=w>0$ on $\partial \Omega$. By the
local maximum principle \cite{m} in $\Omega$ we know that $v+w>0$ in
$\Omega$. This implies that $u>exp(-w)$ in $\Omega$. Hence, we have
the uniform bound
\begin{equation}\label{bound}
exp(-w)<u\leq 1,\ in \  \Omega\times (0,T).
\end{equation}
Then we can extend the solution $u$ beyond any finite time $T>0$.
That is to say, $u$ is a global solution to the equation
(\ref{yamabe}). We denote this solution by $u=u_\Omega$.

We now write $M=\cup_{j=1}^\infty \Omega_j$, $\Omega_j\subset
\subset \Omega_{j+1}$, the exhaustion of bounded smooth domains of
the manifold $M$. Let $u_j=u_{\Omega_j}$. With the help of the
uniform bound (\ref{bound}), we can extract a convergent subsequence
in $C_{loc}^{\infty}(M\times (0,\infty))$ with the limit $u$, which
is the solution to (\ref{yamabe}) on the whole manifold $M$ with the
initial value $u(x,0)=1$.

Again we denote by $g(t)=u(x,t)^{4/(n-2)}g_0$ for the limit solution
$u$. Then by the maximum principle, we know that $R(t)>0$ in
$M\times (0,\infty)$.

In conclusion we have proved the following result.

\begin{Thm}\label{global} Assume that $(M,g_0)$ is a complete non-compact
Riemannian manifold with non-negative scalar curvature $R_0$. Assume
that the Poisson equation $-\Delta w=\frac{n-2}{4(n-1)}R_0$ on $M$
has a non-negative solution $w$. Then the Yamabe flow (\ref{yamabe})
has a global positive solution.
\end{Thm}

From the argument above we can also prove the following result.
\begin{Thm}\label{global+1} Assume that $(M,g_0)$ is a complete non-compact
Riemannian manifold with non-negative scalar curvature $R_0$. Assume
that the Schrodinger equation $-\Delta v+\frac{n-2}{4(n-1)}R_0v=0$
on $M$ has a bounded positive solution $v$. Then the Yamabe flow
(\ref{yamabe}) has a global positive  solution.
\end{Thm}
In fact, we may assume that $v(x)\in (0,1)$. Then, using the
equation (\ref{schrodinger}) to compare  $u$ and $v$, we conclude
that $u(x,t)>v(x)$. Then as above, we can get a global Yamabe flow
$g(t)=u(x,t)^{4/(n-2)}g_0$ such that $v(x)<u(x,t)\leq 1$.

\section{Proof of main result}\label{sect4}
In the first part of this section, we assume that the manifold $M$
satisfies the assumptions in Theorem \ref{main}. Since the initial
metric has bounded curvature, the global Yamabe flow obtained in
previous section has bounded curvature at any finite time interval
(see Theorem 1 in \cite{mcz}). This property makes sure that one can
use Hamilton's tensor maximum principle \cite{H82}.

 To prove
our main result Theorem \ref{main}, we need to state the fundamental
Harnack inequality of B.Chow \cite{chow} for Yamabe flow $g=g(t)$
defined on locally conformally flat manifold $(M,g_0)$. As we
noticed in \cite{macheng}, B.Chow proved this result for the compact
case. However, since one can use Hamilton's tensor maximum
principle, B.Chow's argument can be carried out into complete
non-compact case. Recall that the Harnack quantity defined by B.Chow
is
$$
Z(g,X)=(n-1)\Delta_gR+g(\nabla_gR,X)+\frac{1}{2(n-1)}Rc(X,X)+R^2+\frac{R}{t},
$$
for any metric $g$ and $1-$form $X$. Then the Harnack inequality of
Yamabe flow due to B.Chow is
\begin{Thm}
Let $(M,g_0)$ be a locally conformally flat manifold with positive
Ricci curvature, then under the yamabe flow,
\begin{equation}\label{bchow}
Z(g,X)\geq 0, \ \ t>0,
\end{equation}
for any $1-$form $X$.
\end{Thm}

{\bf Proof of Theorem \ref{main}:} Assume that the scalar curvature
$R_0$ is non-trivial. Using Proposition \ref{ni}, we can solve the
Poisson equation
$$ -\Delta w=\frac{n-2}{4(n-1)}R_0
$$ on $(M,g_0)$ to get the non-negative solution $w$. We then run the Yamabe flow with
the initial metric $g_0$. By Theorem \ref{global}, we have a global
Yamabe flow $g(t)$. As noted in \cite{chow} (see also
\cite{macheng}), the non-negativity property of the Ricci curvature
along the Yamabe flow is preserved. Note that
$$
(n-1)\Delta_gR+R^2=\partial_tR
$$
along the Yamabe flow (see \cite{chow}. The relation (\ref{bchow})
can be written as
$$
\partial_tR+g(\nabla_gR,X)+\frac{1}{2(n-1)}Rc(X,X)+\frac{R}{t}\geq
0.
$$
Choose $X=-d\log R$. We then get
$$
\partial_tR+\frac{R}{t}\geq \frac{|\nabla_gR|^2}{2R}.
$$
Hence, $\partial_t(tR)\geq 0$. Then, for $t\geq 1$ and $\tau\in [
\sqrt{t},t]$,
\begin{equation}\label{monotone}
\tau R(x,\tau)\geq \sqrt{t}R(x,\sqrt{t}).
\end{equation}

Recall that the Yamabe flow can be written as
$\frac{4}{n-2}\partial_t\log u=-R$. Then using the bound
(\ref{bound}) we have for some uniform constant $C(n)>0$,
\begin{equation}\label{bar}\int_0^t Rd\tau=- C(n)\log
u(x,t)\leq C(n)w(x).
\end{equation}
Then using (\ref{monotone}) and (\ref{bar}), we have
\begin{eqnarray*}
&\sqrt{t}R(x,\sqrt{t})\log t
\\
&=\int_{\sqrt{t}}^t\frac{\sqrt{t}R(x,\sqrt{t})}{\tau}d\tau
\\
&\leq \int_{\sqrt{t}}^t\frac{\tau R(x,\tau)}{\tau}d\tau
\\
&\leq \int_0^t R(x,\tau)d\tau
\\
&\leq C(n)w(x).
\end{eqnarray*}
Hence,
$$ 0\leq \sqrt{t}R(x,\sqrt{t})\leq \frac{w(x)}{\log t}.
$$ Then as $t\to\infty$, we have for any $x\in M$, the monotone
non-decreasing quantity has its limit $0$, i.e.,
$$
\sqrt{t}R(x,\sqrt{t})\to 0.
$$
Then we have $R(x,t)=0$ on $M$ at $t=0$. This implies that $(M,g_0)$
is flat. This completes the proof of Theorem \ref{main}.

We now discuss some extensions of Theorem \ref{main}. From the
argument above, we observe that we have actually proved the
following stronger result.

\begin{Thm}\label{main+1} Assume that $(M^n,g_0)$ ($n\geq 3$) is a
locally conformally flat complete non-compact Riemannian manifold
with bounded non-negative Ricci curvature. Assume further that the
Poisson equation $-\Delta w=\frac{n-2}{4(n-1)}R_0$ has a
non-negative smooth solution, where $R_0$ is the scalar curvature of
the metric $g_0$. Then $(M,g_0)$ is flat.
\end{Thm}

By using Theorem \ref{global+1} and the argument above, we can also
prove the following result.

\begin{Thm}\label{main+2} Assume that $(M^n,g_0)$ ($n\geq 3$) is a
locally conformally flat complete non-compact Riemannian manifold
with bounded non-negative Ricci curvature. Assume further that the
Schrodinger equation
$$-\Delta v+\frac{n-2}{4(n-1)}R_0v=0$$ on $M$ has a
positive bounded solution, where $R_0$ is the scalar curvature of
the metric $g_0$. Then $(M,g_0)$ is flat.
\end{Thm}
The detail of proofs of these results are omitted.

\section{Global Yamabe flow}\label{sect5}
 In this section, we shall prove Theorem \ref{mai}.
We first give a few remarks.

\textbf{Remark 1}. It is about the curvature invariant under the
Yamabe flow on a locally conformally flat manifold. At any fixed
point, we can take a normal coordinate system such that at the point
$g_{ij}=\delta_{ij}$ and the Ricci tensor is diagonal
$$ Rc=\lambda_1\bigoplus ...\bigoplus\lambda_n,
$$
with the increasing order $\lambda_1 \leq ...\leq \lambda_n$. Then
the section curvatures (see (2.8) in \cite{chow}) are
$$
K_{ij}=R_{ijji}=\frac{1}{n-2}[\lambda_i+\lambda_j-\frac{R}{n-1}].
$$
where $i\not=j$. From equation (2.13) in \cite{chow} with
$\epsilon=\frac{1}{2(n-1)}$ we get that
$$
\partial_tK_{ij}=(n-1)\Delta K_{ij}+\frac{1}{(n-2)^2}(B_{ii}+B_{jj}),
$$
where
$$
B_{ij}=\mu_1\bigoplus...\bigoplus\mu_n
$$
with
$$
\mu_i=\sum_{k>l,k,l\not=i}(\lambda_k-\lambda_l)^2+(n-2)\sum_{k\not=i}(\lambda_k-\lambda_i)\lambda_i.
$$
From our arrangement, the minimum sectional curvature is $K_{12}$
and in this case
$$
B_{11}+B_{22}=\mu_1+\mu_2,
$$
which is
\begin{eqnarray*}
&\sum_{k>l,k,l\not=1}(\lambda_k-\lambda_l)^2+\sum_{k>l,k,l\not=2}(\lambda_k-\lambda_l)^2+
\\
&(n-2)\sum_{k>2}(\lambda_k-\lambda_1)\lambda_1 + (n-2)(\lambda_2-\lambda_1)\lambda_1\\
&(n-2)\sum_{k>2}(\lambda_k-\lambda_2)\lambda_2
-(n-2)(\lambda_2-\lambda_1)\lambda_2 \\
&
=\sum_{k>l,k,l\not=1}(\lambda_k-\lambda_l)^2+
\\
&(n-2)\sum_{k>2}(\lambda_k-\lambda_1)\lambda_1
+(n-2)\sum_{k>2}(\lambda_k-\lambda_2)\lambda_2 \\
&+\sum_{k>l,k,l\not=2}(\lambda_k-\lambda_l)^2-(n-2)(\lambda_2-\lambda_1)^2
\geq 0,
\end{eqnarray*}
since
$$
\sum_{k>l,k,l\not=2}(\lambda_k-\lambda_l)^2-(n-2)(\lambda_2-\lambda_1)^2\geq
0.
$$
Then we can show by using Hamilton's tensor maximum principle the
following result.

\begin{Pro}\label{rem1} Both the non-negative sectional curvature (or positive
curvature operator)
 condition
 of the locally conformally flat manifold $(M,g_0)$ are preserved
along the LCF Yamabe flow with the initial metric $g_0$.
Furthermore, the positivity property of AVR, and the curvature decay
at infinity of the locally conformally flat manifold $(M,g_0)$ are
also preserved along the Yamabe flow (\ref{yamabe0}).
\end{Pro}
One may compare this result with that of Ricci flow \cite{chow06}.

\textbf{Remark 2}.  In dimension two, the Yamabe flow is the same of
Ricci flow and Theorem \ref{mai} is always true in this case. One
may see \cite{DP} and \cite{M} for more related results. We remark
that in the Ricci flow case, people use the stronger assumption
about the positive curvature operator in the understanding of blow
up limit argument and this assumption plays two roles there. One is
that it is used to keep the trace Harnack inequality hold true for
the Ricci flow, which is true for our Yamabe flow with positive
Ricci curvature on the LCF manifold $M$. The other use of the
positive curvature operator in Ricci flow is to make the metrics
have positive sectional curvature so that one can use Toponogov
comparison theorem in the dimension reduction argument (Theorem 8.46
in \cite{chow06}). In the Yamabe flow on the LCF manifold, the
positive sectional curvature is preserved. Using the same proof of
Theorem 8.46 in \cite{chow06}, we get the dimension reduction result
of the Yamabe flow on the LFC manifold with the replacement of the
assumption of bounded nonnegative curvature operator for Ricci flow
by the assumption of bounded nonnegative sectional curvature for the
Yamabe flow. That is, we have

\begin{Pro}\label{rem2} Let $(M,g(t))$, $t\in(-\infty,\Omega)$, $\Omega>0$
be a complete noncompact locally conformally flat Yamabe flow with
nonnegative sectional curvature. Suppose there exist sequences
$x_i\in M$, $r_i\to\infty$, and $A_i\to\infty$ such that
$d_{g_0}(x_i,p)\geq A_ir_i$ and for all $y\in B_{g_0}(x_i,A_ir_i)$,
$R(y,0)\leq r_i^{-2}$. Assume further that there exists an
injectivity radius lower bound at $(x_i,0)$ in the sense that
$inj(x_i,0)\geq \delta r_i$ for some $\delta>0$. Then a sebsequence
of solutions $(M,r^{-2}g(r_i^2t),x_i)$ converges to a complete limit
solution $(M_\infty, g_\infty(t),x_\infty)$ which is the product of
an $(n-1)$-dimensional solution (with bounded nonnegative sectional
curvature) with a line.
\end{Pro}

\textbf{Remark 3}. We shall use the point picking method of
G.Perelman \cite{p}. For purpose, we recall that for a complete
noncompact Riemannian manifold $(M^n,g)$, the asymptotic scalar
curvature ratio (ASCR, in short) is defined by
$$
ASCR(M,g)=\overline{\lim}_{d(x,0)\to\infty} R(x)d_g(x,0)^2,
$$
where $0\in M$ is any choice of origin. When $M$ is fixed, we write
by $ASCR(g)=ASCR(M,g)$. As noted by Hamilton (p.307 in
\cite{chow06}, the definition does not depend the choice of $0$ and
if $(\tilde{M}, \tilde{g})$ is a Riemannian covering of $(M,g)$,
then $ASCR(\tilde{M},\tilde{g})\geq ASCR(M,g)$. Using the bounded
positive Ricci curvature assumption on the LCF manifold $M$ and
Chow's trace Harnack inequality for the scalar curvature, one can
show (as in Theorem 8.32 in \cite{chow06}) that for the complete
ancient solution to the Yamabe flow $ASCR(g(t))$ is independent of
the time variable $t$ and if $ASCR(g(t))=\infty$, then
$ASCR(g(t_-))=\infty$ for all $t_-\leq t$. Once we have
$ASCR=\infty$ we can have a good understanding of $(M,g(t))$ at
space infinity.

\textbf{Remark 4}. The idea of the proof of Theorem \ref{mai} is
below. We shall argue by contradiction. We shall assume the result
is true in dimension $n-1$ and we shall prove the result in
dimension $n$. If the Yamabe flow only exists up to the finite time
$T>0$, then as in \cite{macheng} we can take a blow-up sequence
$(x_j, t_j)\in M\times (0,T)$, $t_j\to T$ and make a new Yamabe flow
$g_j(x,t)$ with base point $(x_j,t_j)$, which converges in the sense
of Cheeger-Gromov to an ancient Yamabe flow $\bar{g}(t)$,
$-\infty<t<\Omega$ with $0\leq \Omega\leq +\infty$, with the
non-negative sectional curvature (or positive curvature operator)
and the positivity property of AVR. Here we have used the
semi-continuity of AVR under the convergence of Cheeger-Gromov.
However, with the help of differential Harnack inequality of B.Chow
\cite{chow}, arguing as in the work of G.Perelman \cite{p}, we can
show that for the ancient solution $\bar{g}(t)$,
$$
AVR(\bar{g}(t))=0,
$$
which gives us a contradiction with $AVR>0$ for our Yamabe flow.

 Here is the\emph{ proof of Theorem} \ref{mai}.
 \begin{proof}
We shall present the argument without full detail since it is
familiar to experts in Ricci flow and many details can be found in
the book \cite{chow06}(for Ricci flow). Using Hamilton's tensor
maximum principle for the Yamabe flow, we may assume that the
sectional curvature is always positive along the Yamabe flow.
Arguing in the same way as in Theorem 18.2 in \cite{H2}, one know
that the curvature decay at space infinity is preserved, i.e.,
$\lim_{d(x,0)\to\infty}|Rm(g(t))|=0$ for all $0<t<T$. With this
understanding and the fact that our Yamabe flow $(M,g(t))$ has
bounded positive Ricci curvature, one shows (as in Theorem 8.37 in
\cite{chow06}) that $AVR(g(t))>0$ is independent of $t$. Then as in
\cite{macheng} we can get a blow up limit $(\bar{M},\bar{g}(t))$,
$-\infty<t<\Omega, 0<\Omega<\infty$, at the blow-up time $T<\infty$
and with $AVR>0$.

We need consider two cases (see Definitions 4.1 and 4.2 in
\cite{macheng}).

(A). One is that the Yamabe flow $(M,g(t))$ is type I singularity,
i.e, $sup_{M\times (0,T)}(T-t)R<\infty$.

(B). The other case is that Yamabe flow $(M,g(t))$ is type IIa
singularity, i.e, $sup_{M\times (0,T)}(T-t)R=\infty$.

In case (A), the blow up limit $(\bar{M},\bar{g}(t))$,
$-\infty<t<\Omega, 0<\Omega<\infty$, is the ancient Yamabe solution
such that $R(\bar{g}(t))\leq \frac{\Omega}{\Omega-t}$ everywhere
with equality holding somewhere at $t=0$. We claim that
$ASCR(\bar{g}(t))=\infty$ for all $t$. If $ASCR(\bar{g}(t_0))=0$ for
some $t_0$, then we may invoke Theorem B in \cite{PT} to conclude
that the universal cover of $(\bar{M},\bar{g}(t_0))$ is isometric to
$R^{n-2}\times (\Sigma,g_\Sigma(t_0))$ with $ASCR(g_\Sigma(t_0))=0$.
According the strong maximum principle of Hamilton, we know that the
universal cover of the flow $(\bar{M},\bar{g}(t))$ is isometric to
$R^{n-2}\times (\Sigma,g_\Sigma(t))$, where $(\Sigma,g_\Sigma(t))$
is an ancient 2-dimensional Yamabe k-solution (which is also the
Ricci k-solution and k-solution is defined as in \cite{p}). However,
by the result of Hamilton, $(\Sigma,g_\Sigma(t))$ is the round
2-sphere $S^2$ and $ASCR(R^{n-2}\times S^2)$ is not zero, a
contradiction.

If $0<ASCR(\bar{g}(t_0))<\infty$ for some $t_0$, then we have a
positive constant $C>0$ such that
\begin{equation}\label{cheng}
C^{-1}\leq R(x,t_0)(d_{t_0}(x,0)^2+1)\leq C
\end{equation}
on $\bar{M}$, where $R(x,t_0)$ and $d_{t_0}(x,0)$ are the scalar
curvature and distance function of $\bar{g}(t_0)$ respectively. We
new follow the argument in \cite{chow06} (from line -6 in page 347
to line-10 in page 349). Using Chow's Harnack inequality we have
$$
R_t\geq \frac{1}{2} (Rc^{-1})^{ij} \nabla_iR \nabla_jR.
$$
Recall from (9.26) in \cite{chow06} that
$$
Rc^{-1}(\nabla R, \nabla R)\geq \frac{<\nabla R,X>^2}{Rc(X,X)}
$$
for any nonzero vector $X$. Let $\gamma:[0,L]\to \bar{M}$ be a
minimizing geodesic and denote by
$$
\partial_s R=\partial_s R(\gamma(s),t)=<\nabla R,\gamma'(s)>.
$$
Then using the above two inequalities with $X=\gamma'$, we have as
in (9.27) in \cite{chow06} that
\begin{equation}\label{del}
(\partial_s R)^2\leq 2Rc(\gamma',\gamma')R_t.
\end{equation}

Using Chow's trace Harnack inequality again (for the scalar
curvature), $R_t\geq 0$, we know that
$$
R(x,t)(d_{t_0}(x,0)^2+1)\leq C, \ \ on \ \ \bar{M},
$$
for $t\leq t_0$. Since $Rc>0$, we also have
$$
|Rc(x,t)|(d_{t_0}(x,0)^2+1)\leq C, \ \ on \ \ \bar{M}.
$$
Applying Shi's local derivative estimates of the Yamabe flow in the
domain $B_{\bar{g}(t_0)}(x,\frac{1}{2}(1+d_{t_0}(x,0)))$, (see line
-8 in page 348 in \cite{chow06} and the local derivatives are true
for Yamabe flow as noted in \cite{macheng}, see also Theorem 2.2 in
\cite{macheng}), we then know that
$$
\partial_tR(x,t_0)\leq C(d_{t_0}(x,0)^2+1)^{-2}.
$$
Using the first inequality in (\ref{cheng}) we know from
(\ref{scalar}) that
$$
\partial_tR(x,t_0)\leq C R(x,t_0)^2, \ \ on \ \bar{M}.
$$
Thus (\ref{del}) implies that at $t=t_0$,
$$
(\partial_s \log R)^2\leq CRc(\gamma',\gamma').
$$
Then we can use the argument in page 349 in \cite{chow06} to get
that for $1\leq s_1\leq s_2\leq L/2$,
$$
\log (\frac{R(\gamma(s_1))}{R(\gamma(s_2))})\geq c\log^2
(\frac{R(\gamma(s_1))}{R(\gamma(s_2))})-C.
$$
Hence
$$
\frac{R(\gamma(s_1))}{R(\gamma(s_2))}\leq C
$$
for some uniform constant $C>0$. Sending $L\to\infty$ and
$s_2\to\infty$ (and by $ASCR<\infty$, $R(\gamma(s_2))\to 0$), we get
$R(\gamma(s_1))=0$, which is absurd to the fact $R>0$ at
$\gamma(s_1)$. Then we must have
$$
ASCR(\bar{g}(t))=\infty.
$$
In this case we can use the dimension reduction argument
(Proposition \ref{rem2}, similar to Theorem 8.46 in \cite{chow06},
see also case 1 in page 352 in \cite{chow06}) to conclude a
contradiction to the fact that $AVR>0$.

In case (B), we know that the blow up limit is an eternal solution
with property $AVR>0$ and with positive sectional curvature. In this
case, using Chow's Harnack inequality to the quantity
$Z+\frac{R}{t-\alpha}$, where $$ Z:=R_t+<\nabla
R,X>+\frac{1}{2(n-1)}R_{ij}X^iX^j
$$
(and getting $Z\geq 0$ by sending $\alpha\to -\infty$ ) and arguing
as proof of Theorem 4 in \cite{macheng}, we can get that the
backward limit $(\bar{M}_\infty,\bar{g}_\infty(t))$ is a steady
gradient Yamabe soliton with positive sectional curvature. This is
also from Theorem 2.5(i) in \cite{G}. Then the proof of this part is
as in case 2 in page 352 in \cite{chow06} and we only need to
consider $(\bar{M}_\infty,\bar{g}_\infty(t))$ being the steady
soliton. Again, we shall show that $ASCR=\infty$ and we use
dimension reduction as in case (A) to get a contradiction. We denote
by $M$ the manifold in the soliton solution instead of
$\bar{M}_\infty$. Then, there is a smooth concave function $f_0:M\to
R$ and a one parameter diffeomorphism $\phi(t)$ such that
$\bar{g}_\infty(t)=\phi^{*}(t)g_0$ and
\begin{equation}\label{ricci}
R(\bar{g}_\infty)\bar{g}_\infty+Hess_{\bar{g}_\infty}f=0
\end{equation}
where $f=\phi(t)^{*}(f_0)$. Since $R>0$ we know that $f_0$ is
concave and proper. Using the Morse theory, we know that $M$ is
diffeomorphic to $R^n$. Again, as in the proof of Theorem 9.44 in
\cite{chow06}, we can conclude that $ASCR(\bar{g}_0)=\infty$. In
fact, for $0<ASCR(\bar{g}_0)<\infty$, letting $0$ be the maximum
point of $R$ and letting $h=-f_0$, Then $M-\{0\}$ is diffeomorphic
to $S^{n-1}\times R$. Although we don't have the relation $R+|\nabla
f|^2=constant$ as in steady Ricci soliton, we can use the following
argument of Hamilton \cite{H2}(see also the proof of Theorem 1 in
\cite{mv}). Taking the covariant derivative of the relation
(\ref{ricci}) we have
$$
\nabla_iR\bar{g}_{jk}+ \nabla_i\nabla_j\nabla_kf=0.
$$
Since
$$
\nabla_i\nabla_j\nabla_kf=\nabla_j\nabla_i\nabla_kf=R_{ijkl}\nabla_lf,
$$
we get by the contracted Bianchi identity
$$
\nabla_jR_{ij}=\frac{1}{2}\nabla_iR
$$
that $$
\nabla_iR\bar{g}_{jk}-\nabla_jR\bar{g}_{ik}+R_{ijkl}\nabla_lf=0
$$
and
$$
(n-1)\nabla_i R=2R_{ij}\nabla_jf.
$$
Hence $\nabla f(0)=0$ and $0$ is the unique maximum point of $f$.

 We now take a minimizing geodesic $\gamma(s)\subset (M,g_0)$ starting from $0$ and let $X=\gamma'$. For any $\epsilon>0$,
  we can use the relation
$$
X (X h)=D^2h(X,X)=R|X|^2=R>0
$$
to conclude that
$$
X(h)(s)=\int_0^s R(\gamma)\geq C_0(\epsilon)>0
$$
for some uniform constant $C_0=C_0(\epsilon)>0$, depending only on
$\inf_{B_{g_0}(0,\epsilon)} R$, and for all $s\geq \epsilon$. Then
$|\nabla h|\geq C_0>0$. By this and arguing as in page 355 in
cite{chow06} we can get the key estimate
$$
d_{g_0}(\phi(t),0)\geq c|t|
$$
for all $x\in M-B_{g_0}(0,\epsilon)$ as in (9.33) in \cite{chow06}.
We can use the argument in page 355 to page 356 in \cite{chow06} to
get a contradiction. Hence we have $ASCR=\infty$ and we use the
dimension reduction as in case (A) to get a contradiction again.
Therefore $T=\infty$ and we have the global Yamabe flow $g(t)$.

This completes the proof of Theorem \ref{mai}. \end{proof}

We may give some remarks about the proof above. From the proof
above, we actually prove the following results.

\begin{Thm}\label{mai1}
Assume that $(M,g(t))$, $t\in (-\infty,0]$ is a locally conformally
flat type I ancient solution to Yamabe flow with bounded positive
sectional curvature and $ASCR(g(t))<\infty$. Then $AVR(g(t))>0$ and
for any choice of  origin $0\in M$, there exists a constant
$c=c(0,t)>0$ such that $$ R(x,t)(1+d_{g(t)}(x,0))^2\geq c, \ \ on \
M.$$
\end{Thm}
This result is similar to Proposition 9.22 in \cite{chow06}.
\begin{Thm}\label{mai2}
Assume that $(M,g(t))$, $t\in (-\infty,0]$ is a locally conformally
flat ancient solution to Yamabe flow with bounded nonnegative
sectional curvature. Then $AVR(g(t))=0$.
\end{Thm}
This result is similar to Proposition 9.30 in \cite{chow06}.
\begin{Thm}\label{mai3}
Assume that $(M,g(t))$, $t\in (-\infty,\Omega)$ is a locally
conformally flat type I ancient solution to Yamabe flow with bounded
nonnegative sectional curvature. Then $ASCR(g(t))=\infty$ for all
$t$.
\end{Thm}
This result is similar to Proposition 9.32 in \cite{chow06}.

We conclude the paper by mentioning a question. It is interesting to
consider the behavior of the global Yamabe flow as $t\to\infty$. The
dimension two case is well done in some cases \cite{M}.

\emph{Acknowledgement}:  The author would like to thank the unknown
referees very much for helpful suggestions.


\begin{thebibliography}{20}
\bibitem{AM}
 Y.An, L.Ma, \emph{The Maximum Principle and the Yamabe
Flow}, Partial Differential Equations and Their Applications, World
Scientific, Singapore, pp211-224, 1999


\bibitem{chow} B.Chow,
\emph{Yamabe flow on locally conformally flat manifolds}, Comm. pure
appl. math., Vol.XLV,(1992)1003-1014


\bibitem{chow06} B Chow, P. Lu, L. Ni, \emph{Hamilton's Ricci flow}, Graduate Studies in
Math. 77, Amer. Math. Soc. (2006) MR2274812

\bibitem{DP}
P.Daskalopoulos, M. del Pino, \emph{Type II collapsing of maximal
solutions to the Ricci flow in R2}. Ann. Inst. H. Poincare Anal. Non
Lineaire 24 (2007), no. 6, 851-874.

\bibitem{DS}
P.Daskalopoulos, N.Sesum,  \emph{The classification of locally
conformally flat Yamabe solitons}, arxiv.1104.2242v1

\bibitem{G} Gu, H.L.: \emph{Manifolds with pointwise Ricci pinched curvature}. Acta
Math.Scientia, 30B(3)(2010)819-829

\bibitem{H82}
R.Hamilton, \emph{Three-manifolds with positive Ricci curvature}.
 J. Differential Geom., 2(1982)255-306

  \bibitem{H2}
R. S. Hamilton, \emph{Formation of singularities in the Ricci flow},
Surveys in Diff. Geom. 2 (1995), 7-136.

\bibitem{ly}P.Li, S.T.Yau,
\emph{On the parabolic kernel of the Schrodinger operator}, Acta
Math. 156(3-4), 153-201(1986)

\bibitem{ni}L.Ni,
\emph{An optimal gap theorem}, Invent. Math.,Volume 189,
737-761(2012)

\bibitem{M} Li Ma, \emph{Convergence of Ricci flow on R2 to the plane}, Differential
Geometry and its Applications 31 (2013) 388-392

\bibitem{ma}
L. Ma, \emph{Expanding Ricci solitons with pinched Ricci curvature},
KODAI Math. J. 34 (2011), 140-143

\bibitem{macheng} L.Ma, L.Cheng,\emph{
Yamabe flow and Myers type theorem on complete manifolds},
J.Geom.Anal.,DOI 10.1007/s12220-012-9336-y, online, 2012

\bibitem{mcz}
L.Ma, L.Cheng, A.Zhu, \emph{Extending Yamabe flow on complete
Riemannian manifolds}, Bull. Sci. math. 136(2012)882-891

\bibitem{mv}
Li Ma, M. Vicente, \emph{Remarks on scalar curvature of Yamabe
solitons}, Ann Glob Anal Geom., (2012) 42:195-205

\bibitem{m}
C.B.Morrey, Jr., \emph{Multiple Integrals in Calculus of
Variations}. Springer, New York (1966)

\bibitem{p} G.Perelman, \emph{The entropy formula for the Ricci
flow and its geometric applications}. arXiv:math.DG/0211159.

\bibitem{SY}R.Schoen, S.T.Yau, \emph{Conformally flat manifolds, Kleinian groups
and scalar curvature}, Invent. math. 92, 47-71(1988)


\bibitem{shi}W.X.Shi, \emph{Ricci flow and the uniformization on complete
noncompact K\"{a}ler manifolds}, J.Diff.Geom., 45 (1997)94-220

\bibitem{PT}
A.Petrunin, W.Tuschmann, \emph{Asymptotic flatness and the cone
structure at infinity}, Math Ann. 321(2001)775-788.

\bibitem{V} N. Varopoulos, \emph{The Poisson kernel on positively curved
manifolds}, J. Functional Analysis 44(1981)359-380

\end{thebibliography}
\end{document}